\documentclass[11pt,leqno]{amsart}

\usepackage{amsthm,amsfonts,amssymb,amsmath,oldgerm,mathrsfs,mathabx,enumitem,mathtools}
\usepackage[scr=boondoxo]{mathalfa}
\numberwithin{equation}{section}
\usepackage{fullpage,setspace,fancyhdr,bm}
\usepackage{doi}
\usepackage{graphicx,psfrag,listings} 
\usepackage{color}

\usepackage[top=30mm,bottom=30mm,left=25mm,right=25mm,a4paper]{geometry}

\usepackage{hyperref}
\usepackage{graphics}

\setlist[itemize,1]{label=\ensuremath{\diamond}} 
\setlist[enumerate,1]{label=(\textit{\roman*})}

\theoremstyle{plain}
\newtheorem{theorem}{Theorem}

\newtheorem{proposition}[theorem]{Proposition}
\newtheorem{lemma}[theorem]{Lemma}

\theoremstyle{definition}

\newtheorem{remark}[theorem]{Remark}

\newtheorem{assumption}[theorem]{Assumption}

\newcommand{\cai}[1]{\mathbf{1}_{#1}}

\newcommand{\buc}[1]{\ensuremath{BUC^{#1}}}
\newcommand{\ad}{\mathrm{ad}}

\newcommand\dd{\,\mathrm{d}}

\newcommand{\deq}{:=}
\newcommand{\re}{\mathrm{e}}
\newcommand{\ri}{\mathrm{i}}

\DeclareMathOperator{\Real}{Re}
\DeclareMathOperator{\Imag}{Im}

\newcommand{\smallo}{\operatornamewithlimits{o}}

\newcommand{\norm}[1]{\lVert#1\rVert}
\newcommand{\Norm}[1]{\left\lVert#1\right\rVert}
\newcommand{\absolute}[1]{\lvert#1\rvert}
\newcommand{\Absolute}[1]{\left\lvert#1\right\rvert}
\newcommand{\scalp}[1]{\langle #1\rangle} 
\newcommand{\Scalp}[1]{\left\langle #1\right\rangle}

\newcommand*{\Set}[1]{\left\lbrace #1 \right\rbrace}

\newcommand{\CC}{{\mathbb C}}

\newcommand{\NN}{{\mathbb N}}

\newcommand{\RR}{{\mathbb R}}

\newcommand{\ZZ}{{\mathbb Z}}

\newcommand\uphi{{\underline \phi}}

\newcommand\cA{{\mathcal A}}

\newcommand\cC{{\mathcal C}}

\newcommand\cL{{\mathcal L}}

\newcommand\cO{{\mathcal O}}

\newcommand\cQ{{\mathcal Q}}
\newcommand\cR{{\mathcal R}}

\title{Existence of monostable fronts for a KPP infinite-difference numerical scheme}

\author{Louis Garénaux}
\address{Englerstraße 2, 76131 Karlsruhe, Germany}
\email{{\tt louis.garenaux@kit.edu}}
\thanks{LG warmly thanks the Leiden Mathematisch Instituut for providing everything necessary during its stay in Leiden. The research stay was funded by Karlsruhe House of Young Scientist (KHYS). Research of LG was Funded by the Deutsche Forschungsgemeinschaft (DFG, German Research Foundation) – Project-ID 258734477 – SFB 1173}

\author{Hermen Jan Hupkes}
\address{Niels Bohrweg 1, 2333 CA Leiden, The Netherlands}
\email{{\tt hhupkes@math.leidenuniv.nl}}

\begin{document}

\begin{abstract}
We study the existence of traveling wave solutions for a numerical counterpart of the KPP equation. We obtain the existence of monostable fronts for all super-critical speeds in the regime where the spatial step size is small. The key strategy is to transfer the invertibility of certain linear operators related to the front solutions from the continuous setting to the discrete case we are interested in. We rely on resolvent bounds which are uniform with respect to the step size, a procedure which is also known as spectral convergence. The approach is also able to handle infinite range discretizations with geometrically decaying coefficients that are allowed to have both signs, which prevents the use of the comparison principle.

\vspace{0.5em}


\vspace{0.5em}

{\small \paragraph {\bf AMS Subject Classifications (MSC 2020):} 39A14, 39A12, 39A70, 47A10, 47N40, 47A12.}
\end{abstract}

\maketitle

\section{Introduction}

\subsection{Problem and statements}

Our primary goal is to construct traveling front solutions 
to the equation
\begin{equation}
\label{e:main}
\dot{u}_j(t) = \frac{u_{j+1} - 2 u_j + u_{j-1}}{h^2} + g(u_j), \hspace{4em} t>0, \quad j\in \ZZ.
\end{equation}
Equation \eqref{e:main} appears as a space discretization of a scalar reaction-diffusion equation: we restrict to small step size $h>0$. More precisely, we focus on monostable equations, by which we mean in this work that $g:\RR \to \RR$ satisfies for an integer $K\geq 2$ the 
\begin{assumption}
\label{a:g}
The function $g\in \cC^K$ admits two consecutive zeros $u=0$ and $u=1$ and is concave on $[0,1]$. It is non-degenerate, in the sense that $g'(1) < 0 < g'(0)$.
\end{assumption} 

The solution we are aiming at is a traveling heteroclinic connection between these two states, meaning that we look for a solution 
\begin{equation}
\label{e:traveling-wave}
u_j(t) = \phi_h(jh - ct), \hspace{4em} t>0, \quad j\in \ZZ,
\end{equation}
defined by its speed $c>0$, and its profile $\phi_h \in \cC^{K+1}(\RR,\RR)$ that satisfies
\begin{equation}
\label{e:heteroclinic-connection}
\lim_{-\infty} \phi_h = 1, \hspace{4em}
\lim_{+\infty} \phi_h = 0.
\end{equation}

\begin{theorem}
\label{t:local}
Assume that Assumption \ref{a:g} holds, and let $c>2\sqrt{g'(0)}$. Then there exists $h_0>0$ such that for all $h \in (0, h_0)$, equation \eqref{e:main} admits a solution of the form \eqref{e:traveling-wave}-\eqref{e:heteroclinic-connection}.
\end{theorem}

The previous statement partially recovers a result actually known since \cite{Zinner-Harris-Hudson-93}. Their work sparked a number of follow-up studies focusing on the existence, monotonicity, uniqueness and stability of monostable front solutions to semi-discrete difference equations \cite{Carr-Chmaj-04,Ma-Zou-05}. In addition, various extensions have been considered involving infinite difference schemes \cite{Ma-Weng-Zou-06}, time dependent coefficients \cite{Ducrot-Jin-24}, non-linear porous-media type diffusions \cite{Chen-Guo-02} and interaction terms $g(u)$ featuring time-delays and spatial convolutions \cite{Ma-Zou-05}. All these results make critical use of the comparison principle satisfied by \eqref{e:main}.

The techniques we use in the present paper are fundamentally different in nature and are based on the so-called `spectral convergence' approach developed in \cite{Bates-Chen-Chmaj-03}. Studying bistable equations such as \eqref{e:main} with $g(u) = u(1-u)(u-a)$ and $a \in (0,1)$, they managed to transfer the existence result of the continuous problem ($h\to 0$) to the discrete problem $(h>0)$. This motivated several extensions \cite{Hupkes-VanVleck-22, Hupkes-Jukic-24}, showcasing that results concerning the existence and stability of waves and the asymptotic description of their properties can be obtained for a wide range of bistable settings. In particular, infinite difference schemes \cite{Schouten-Straatman-Hupkes-21}, equations posed on trees \cite{Hupkes-Jukic-24}, reaction-diffusion systems of FitzHugh-Nagumo type \cite{Schouten-Straatman-Hupkes-19} and time dependent non-uniform discretization grids \cite{Hupkes-VanVleck-23} were studied.

We see the above theorem as an illustration that the approach of the second branch of literature
can be applied to the monostable problems studied in the first branch. As a direct consequence, we are able to deal with a larger class of discretization schemes in the monostable setting. In particular, let us consider an infinite difference approximation of the second order derivative:
\begin{equation}
\label{e:main-2}
\dot{u}_j(t) = \sum_{k=1}^{+\infty} a_k \frac{u_{j+k} - 2 u_j + u_{j-k}}{k^2 h^2} + g(u_j), \hspace{4em} t>0, \quad j\in \ZZ,
\end{equation}
in which the sequence $(a_k)_{k \geq 1} \subset \RR$ satisfies the 
\begin{assumption}
\label{a:a}
The sequence is geometrically decaying: there exist $C>0$ and $\rho\in (0,1)$ such that for all $k\geq 1$, $\absolute{a_k} \leq C\rho^k$. The sum of the coefficients is positive, and normalized as $\sum_{k=1}^{+\infty} a_k = 1$.
\end{assumption}
Observe that equation \eqref{e:main} corresponds to the choice $a_1 = 1$ and $a_k = 0$ for all $k\geq 2$. However, our techniques do not require the hypothesis that $a_k \geq 0$ for all $k\geq 1$, which is necessary for the comparison principle to hold. In particular, our main result here goes beyond previous infinite-range results \cite{Ma-Weng-Zou-06} that were established using the comparison principle. 

\begin{theorem}
\label{t:inf-range}
Assume that both assumptions \ref{a:g} and \ref{a:a} hold, and let $c>2\sqrt{g'(0)}$. Then there exists $h_0>0$ such that for all $h\in (0, h_0)$, equation \eqref{e:main-2} admits a solution of the form  \eqref{e:traveling-wave}-\eqref{e:heteroclinic-connection}.
\end{theorem}

We remark that discrete convolution kernels with negative coefficients arise naturally in applications, such as particle chains with visco-elastic interactions, see \cite[Eq. (11)]{Vainchtein-VanVleck-09}. Furthermore, the use of non-positive $a_k$ allows one to consider numerical differentiation schemes that have higher orders of consistency. For example, $(a_k)_{k \geq 1} = (-\frac{1}{3},\frac{4}{3},0,0, \dots)$ satisfies Assumption \ref{a:a}, and corresponds to the discrete diffusion operator
\begin{equation*}
\frac{-u_{j+2} + 16 u_{j+1} - 30 u_j + 16 u_{j-1} - u_{j-2}}{12 h^2} \simeq u''(jh) + \cO(h^4),
\end{equation*}
which should be contrasted to the $\cO(h^2)$-approximation made in \eqref{e:main}. Another interesting example is given by $(a_k)_{k \geq 1} = (-\frac{1}{2},\frac{3}{2},0,0, \dots)$, since it corresponds to a convolution kernel $\frac{1}{h^2}(\frac{3}{8}, -\frac{1}{2}, \frac{1}{4}, -\frac{1}{2}, \frac{3}{8})$ that has a positive central coefficient $\frac{1}{4h^2}$.

\subsection{Obstacles and approach}

The main idea of our proof is to rely on the corresponding spatially continuous equation
\begin{equation}
\label{e:continuous-pb}
\partial_t u = \partial_{xx} u + g(u), \hspace{4em} t>0, \quad x\in \RR,
\end{equation}
for which existence of traveling waves $u(t,x) = \phi_0(x - ct)$ satisfying \eqref{e:heteroclinic-connection} is well-known \cite{Kolmogoroff-Pretrovsky-Piscounoff-37,Fisher-37} for all speeds $c\geq 2\sqrt{g'(0)}$ and under a slightly weaker version of assumption \ref{a:g}. We construct a fix-point argument for the space-discretized equation, linearizing near an approximate solution. When $h$ is small and $c\neq 2\sqrt{g'(0)}$, we are able to transfer the invertibility properties of the linear part from the space-continuous setting towards the space-discretized setting.
Throughout the paper we only consider the infinite-range setting of Theorem \ref{t:inf-range}, but for the remainder of our discussion here we will focus on the setting of Theorem \ref{t:local}.

As we will see, our search for discrete traveling front solutions replaces the time-derivative $\partial_t$ in \eqref{e:continuous-pb} by the convection operator $c\partial_x$, while the Laplacian $\partial_{xx}$ is replaced by the discrete version $\Delta_{kh}$ that acts as 
\begin{equation*}
\Delta_{kh} v \deq \frac{1}{k^2 h^2} \bigg(v(\cdot + kh) - 2 v + v(\cdot - kh)\bigg),
\end{equation*}
In particular, we can view the combination 
$\Delta_{kh} + c\partial_x$ as an $L^2(\RR)$-based operator with domain $H^1(\RR)$. In this setting, the limit $h\to 0$ is singular. To obtain the invertibility of the linear operator in the discrete setting, we establish a $h$-uniform resolvent bound, by relying on the $L^2$-weak convergence of $\Delta_{kh} v$ towards $v''$. 

In the present monostable setting, the approach initiated in \cite{Bates-Chen-Chmaj-03} for bistable waves does not immediately carry over, and two main points have to be addressed.

\begin{itemize}
\item The bistable structure localizes solutions of the resolvent problem by preventing the bulk of their energy from escaping to infinity as $h \to 0$. In the monostable setting, a similar structure can be recovered using over-localizing weights. However, these introduce delicate novel terms when interacting with the shift operators, which need to be carefully controlled. The geometric decay of the sequence $(a_k)_k$ is crucial to achieve this control; see Proposition \ref{p:Delta-h-decomposition} below.

\item The residual terms that appear when linearizing the profile equation close to an approximate solution are not localized enough to use the aforementioned weights. We use a tailor-made approximate solution to cancel the most critical residual terms, following the so-called far-field approach in \cite{Avery-Garenaux-23}.
\end{itemize}

Our inspiration to introduce spatial weight functions stems from the continuous problem \eqref{e:continuous-pb}, where such weights were used by Sattinger in \cite{Sattinger-76} to counteract the instability of the equilibrium $u = 0$. However, they are  available for fast waves $c>2\sqrt{g'(0)}$ only. For the slowest wave $c=2\sqrt{g'(0)}$, no spectral gap can be recovered \cite{Faye-Holzer-19} since a branch point of the absolute spectrum lies on the imaginary axis \cite{Sandstede-Scheel-00,Faye-Holzer-Scheel-Siemer-22}.

\subsection{Future directions}
We view this paper as a first proof of concept that the spectral convergence method can be applied in monostable settings. As a consequence, it naturally opens several follow up questions.
\begin{itemize}
\item The main obstacle towards understanding the slowest wave $c=2\sqrt{g'(0)}$ is the possibility for the mass to escape at $+\infty$, since exponential weights can no longer be used to stabilize the spectrum. We believe that this can be resolved by a careful analysis of the asymptotic behavior of solutions to the  advance-delay differential equations underlying the resolvent problems for $h > 0$.

\item Most efficient modern solvers for parabolic systems feature time-dependent numerical schemes, where the density of grid points is increased in areas where the solution is steep. Naturally, it is highly relevant to understand the impact that such adaptive schemes have on the properties of the underlying solution. In the bistable case, the persistence of traveling waves has been established in the series of papers \cite{Hupkes-VanVleck-22a,Hupkes-VanVleck-22,Hupkes-VanVleck-23}, which rely heavily upon the spectral convergence approach.

\item In the context of numerical analysis, it is of fundamental importance to study the behaviour of fully discretized schemes such as 
\begin{equation*}
\frac{u_j^{n+1} - u_j^n}{h_t} = \frac{u_{j+1}^{n+1} - 2 u_j^{n+1} + u_{j-1}^{n+1}}{{h_x}^2} + g(u_j^n),
\end{equation*}
which arises from \eqref{e:main} by applying the implicit Euler scheme to the temporal variable. In particular, we are interested in the impact of such schemes on the existence, uniqueness and stability of traveling fronts. Partial answers in several bistable settings have been obtained in \cite{hupkes2016travelling,Schouten-Straatman-Hupkes-21}, again using the spirit of the spectral convergence approach.
\end{itemize}

\section{Notations}

\subsection{Spaces}
In the following, we simply write $L^p$ to denote the usual Lebesgue vector spaces $L^p(\RR)$. When using other spatial domains than $\RR$, we will be more explicit.

\subsection{Taylor remainders}
For any $j\in \NN$ and $f\in \cC^j(\RR, \RR)$, we abbreviate the Taylor expansion at $a$ in the direction $b$ as
\begin{equation}
\label{e:Taylor}
T_{j,f}(a,b) \deq f(a + b) - \sum_{k=0}^{j-1} \frac{f^{(k)}(a)}{k!} b^k.
\end{equation}
We can already use this notation to check that the numerical scheme we consider is \emph{consistent}.
\begin{lemma}
\label{l:Taylor-expansion}
Assume that $f\in \buc{3}(\RR, \RR)$, and denote
\begin{equation}
\label{e:def-Delta-a-h}
\Delta_{a,h} f \deq \sum_{k=1}^{+\infty} a_k \frac{f(\cdot + kh) - 2 f + f(\cdot - kh)}{k^2 h^2}.
\end{equation}
Then the following statements hold.
\begin{enumerate}
\item \label{i:taylor-1} There exists a positive constant $C$ such that for all $h\in (0,1)$,
\begin{equation*}
\norm{\Delta_{a,h} f - f''}_{L^\infty} \leq Ch .
\end{equation*}
\item \label{i:taylor-2} Assume furthermore that $f$ has exponential behavior at infinity, in the sense that there exist positive constants $h_0$, $C$ and $\kappa$ such that for all $x\in \RR$, all $k \geq 1$ and all $h \in (0, h_0)$
\begin{equation*}
\sup_{\absolute{s}\leq kh} \Absolute{\frac{f^{(3)}(x + s)}{e^{-\kappa x}}} \leq C^{1+kh}.
\end{equation*}
Then there exist positive constants $\tilde{h}_0$ and $C$ such that for all $h\in (0,\tilde{h}_0)$,
\begin{equation*}
\Norm{\frac{\Delta_{a,h} f - f''}{e^{-\kappa \cdot}}}_{L^\infty} \leq C h .
\end{equation*}
\end{enumerate}
\end{lemma}
\begin{proof}
It is direct to check that 
\begin{equation*}
\frac{f(\cdot + kh) - 2 f + f(\cdot - kh)}{k^2 h^2} = f'' + \frac{1}{k^2 h^2} \bigg(T_{3,f}(\cdot, kh) + T_{3,f}(\cdot, -kh)\bigg).
\end{equation*}
Using the integral formulation of the Taylor remainder provides the bound $\absolute{T_{3,f}(x,\pm kh)} \leq \frac{1}{6} k^3 h^3 \norm{f'''}_{L^\infty}$ for all $x\in \RR$. This proves the first point. The second point is almost identical. By taking $h$ sufficiently small, we can enforce the geometric decay
\begin{equation*}
\Absolute{a_k \frac{1}{k^2 h^2} \frac{T_{3,f}(\cdot, kh)}{e^{-\kappa \cdot}}} \leq  k h \absolute{a_k}
\sup_{s\in (x, x+kh)} \Absolute{\frac{f^{(3)}(s)}{e^{-\kappa x}}} \leq C \, k h \, (\rho C^h)^k,
\end{equation*}
which allows us to sum over $k \geq 1$.
\end{proof}

\subsection{Far field approach}
To describe the behavior at infinity, we introduce a partition of unity. Let $\cai{-}$ and $\cai{+}$ be two smooth functions $\RR\to \RR$, that are constant outside of a bounded set:
\begin{equation*}
\cai{-}(x) = 
\begin{cases}
1 & x\leq -1\\
0 & x\geq 1
\end{cases},
\hspace{6em}
\cai{+}(x) = 1 - \mathbf{1}_-(x) .
\end{equation*}

\begin{lemma}
\label{l:decomp-profile}
Let $c>2\sqrt{g'(0)}$. There exists positive constants $\kappa_0$, $\delta$, and $w_0\in L^2(\RR)$, such that the wave profile $\phi_0$ for the continuous problem \eqref{e:continuous-pb} can be decomposed as 
\begin{equation*} 
\phi_0(x) = \mathbf{1}_-(x) + w_0(x) + \mathbf{1}_+(x) \, \phi_0(0) \, \re^{-\kappa_0 x},
\end{equation*}	
in which $w_0$ satisfies the decay rates
\begin{equation*}
\norm{x\mapsto \re^{-\delta x} w_0(x)}_{W^{2,\infty}(-\infty, 0)} + \norm{x\mapsto \re^{(\kappa_0 + \delta) x} w_0(x)}_{W^{2,\infty}(0,+\infty)} < \infty.
\end{equation*}
\end{lemma}
\begin{proof}
Since both $\lim_{\pm\infty} \phi_0$ are hyperbolic equilibria for the ODE
\begin{equation*}
-c\phi' = \phi'' + g(\phi),
\end{equation*}
$\phi_0$ decays exponentially at $\pm\infty$, see \cite[\S 6]{Sattinger-76}. We denote the spatial decay rate at $+\infty$ by $\kappa_0$. The claimed decomposition uniquely defines $w_0$, and the bounds are obtained by choosing $\delta$ small enough and exploiting the fact that $\kappa_0$ is a simple spatial eigenvalue. 
\end{proof}

\subsection{Over-localized weights}
To invert the linear operator, we will first need to remove the essential spectrum at the origin. As is classical for monostable waves, we do so using exponential weights.

\begin{lemma}
\label{l:construct-theta}
Let $c>2\sqrt{g'(0)}$ and recall the constants $\delta$ and $\kappa$ from Lemma \ref{l:decomp-profile}. Then there exists $\theta\in (\kappa_0, \kappa_0 + \frac{\delta}{2})$ such that 
\begin{equation}
\label{e:stable-spectrum}
\theta^2 - c\theta + g'(0) < 0,
\end{equation}
and
\begin{equation*}
c - 2\theta < 0.
\end{equation*}
\end{lemma}
\begin{proof}
The spatial decay rate $\kappa_0$ is given by
the smallest solution of the polynomial equation
\begin{equation*}
\kappa_0^2 - c\kappa_0 + g'(0) = 0;
\end{equation*}
see again \cite[\S 6]{Sattinger-76}.
In particular, solutions to the inequality \eqref{e:stable-spectrum} are all larger than $\kappa_0$. Explicitely computing that 
$\kappa_0 \deq \frac{c}{2} - \sqrt{\frac{c^2}{4} - 4 g'(0)}$, 
we are able to choose $\theta$ sufficiently close to $\kappa_0$ in order to ensure that 
\begin{equation*}
\theta < \min\left(\frac{c}{2}, \kappa_0 + \frac{\delta}{2}\right),
\end{equation*} 
which concludes the proof.
\end{proof}

We fix such a $\theta$ for the rest of the article, and as a first step we introduce the initial  weight function
\begin{equation}
\label{e:weight}
\tilde{\omega}(x) = 
\begin{cases}
1 & x\leq -1\\
\exp\big(-\theta \frac{(x+1)^2}{4}\big) & x\in(-1,1)\\
\exp(-\theta x) & x\geq 1 .
\end{cases}.
\end{equation}
Although this expression already enjoys several convenient properties, it lacks smoothness. Thus, we further refine it with the following statement.
\begin{lemma}
\label{l:construct-omega}
Let $\varepsilon>0$. There exists $\omega \in \cC^\infty(\RR, (0,1])$ such that $\omega(x) = \tilde{\omega}(x)$ when $x\notin (-2,2)$, and such that $\Norm{\omega - \tilde{\omega}}_{W^{2,\infty}(\RR)} \leq \varepsilon$.
\end{lemma} 
\begin{proof}
We set out to smoothen the function $\tilde{\omega}$ on the interval $[-2, 2]$. First, we set a partition of unity by choosing $\chi\in \cC^\infty(\RR, [0,1])$ in such a way that 
\begin{equation*}
\chi(x) = 
\begin{cases}
1 & \text{ when } x \in \left(-\frac{3}{2},\frac{3}{2}\right),\\
0 & \text{ when } x \notin (-2,2).
\end{cases}
\end{equation*}
Then, we approximate $\tilde{\omega}_{|[-2,2]}$ using the density of smooth functions: for $\varepsilon>0$, there exist $\varpi\in \cC^\infty([-2,2], \RR)$ such that
\begin{equation*}
\norm{\tilde{\omega} - \varpi}_{W^{2,\infty}(-2,2)} \leq \varepsilon.
\end{equation*}
Finally, we define $\omega$ as the barycenter
\begin{equation*}
\omega = \chi \varpi + (1-\chi) \tilde{\omega}.
\end{equation*}
It satisfies the claimed properties, and the proof is complete.
\end{proof}

For the remainder of the article, we choose $\varepsilon$, and write $\omega$ for the associated weight. We may decrease $\varepsilon$ and continue using the $\omega$ notation without explicitly stating the $\varepsilon$ dependence.

\subsection{Residual, linear operator and quadratic terms}
Inserting the traveling wave ansatz \eqref{e:traveling-wave} into \eqref{e:main-2}, we see that $\phi_h$ must satisfy
\begin{equation}
\label{e:profile}
-c\phi' = \Delta_{a,h} \phi + g(\phi),
\end{equation}
with $\Delta_{a,h}$ defined in \eqref{e:def-Delta-a-h}.

We start by prescribing the behavior of the profile at infinity, in the spirit of Lemma \ref{l:decomp-profile}. 
In particular, we introduce the function
\begin{equation*}
\phi^\infty(x) = \mathbf{1}_{-}(x) + w_0 + \mathbf{1}_{+}(x) \, \phi_0(0) \, \re^{-\kappa_h x} ,
\end{equation*}
in which the exponent $\kappa_h > 0$ will be chosen later.

We further introduce the notation
\begin{equation*}
\cA_h(\uphi) \deq \Delta_{a,h} + c\partial_x + g'(\uphi)
\end{equation*}
to denote the linearization of \eqref{e:profile} at an approximate solution $\uphi$.

\begin{lemma}
\label{l:v-formulation}
Consider a solution to \eqref{e:profile} that can be written in the form
\begin{equation}
\label{e:perturbed-ansatz}
\phi_h = \phi^\infty + \omega v.
\end{equation}
Then $v$ satisfies 
\begin{equation*}
0 = \cR + \cL_h v + \cQ(v),
\end{equation*}
with
\begin{align*}
\cR \deq {} & \omega^{-1} \cA_h(0)\phi^\infty + \omega^{-1} T_{2,g} (0, \phi^\infty),\\
\cL_hv \deq {} & \omega^{-1} \cA_h(\phi^\infty)(\omega v),\\
\cQ(v) \deq {} & \omega^{-1} T_{2,g}(\phi^\infty, \omega v).
\end{align*}
\end{lemma}
\begin{proof}
We first rewrite \eqref{e:profile} as 
\begin{equation*}
0 = \cA_h(0) \phi + g(\phi) - g'(0)\phi,
\end{equation*}
and insert the decomposition \eqref{e:perturbed-ansatz}. 
We force the apparition of a Taylor remainder of order two, and group together linear terms to arrive at
\begin{equation*}
0 = \cA_h(0) \phi^\infty + T_{2,g}(\phi^\infty, \omega v) + \bigg(\cA_h(0) + g'(\phi^\infty) - g'(0)\bigg) \omega v + \bigg(g(\phi^\infty) - g'(0) \phi^\infty\bigg).
\end{equation*}
Artificially inserting a $g(0) = 0$ inside the last pair of parentheses, we obtain another Taylor remainder of order two. Dividing by $\omega$ concludes the proof.
\end{proof}

\subsection{Choice of ansatz}
A quick inspection of $\cR$ reveals that its second term behaves as $x\mapsto \re^{(\theta-2\kappa_h)x}$ when $x\to +\infty$. On the other hand, its first term (typically) behaves as $x\mapsto \re^{(\theta-\kappa_h)x}$.\footnote{This latter fact is easy to see when dealing with finite differences, \emph{e.g.} for \eqref{e:main}. We refer to the precise computations that follow for the general case.} With the intuition that $\kappa_h$ will represent the spatial decay rate of $\phi_h$ and thus will get close to $\kappa_0$ in the $h\to 0$ limit, we see that the first term in $\cR$ \emph{a priori} has the potential to grow as $\xi \to \infty$, preventing the use of over-localizing weights. We choose $\kappa_h$ so that this problematic asymptotic term vanishes identically, by looking for zeros of the function 
\begin{equation*}
G(h, \kappa) \deq g'(0) - c\kappa  + \sum_{k=1}^{+\infty} a_k \frac{\re^{-\kappa k h} - 2 + \re^{\kappa k h}}{k^2 h^2}.
\end{equation*}

\begin{lemma}
\label{l:definition-kappa}
There exists $h_0>0$ and a map $\kappa:(0,h_0) \to \RR$ such that for any $h \in (0, h_0)$, $G(h, \kappa(h)) = 0$. For this choice of spatial exponential rate, we have
\begin{equation*}
\cA_h(0) (x\mapsto \re^{-\kappa(h) x}) = 0.
\end{equation*}
\end{lemma}
\begin{proof}
Thanks to the geometric decay of the sequence $a_k$, we can compute that 
\begin{equation*}
\cA_h(0) (\re^{-\kappa \cdot}) (x) = \re^{-\kappa x} G(h, \kappa)
\end{equation*}
holds for all $x\in \RR$, establishing the final statement. We now focus on the definition of $\kappa(h)$. Using Taylor expansions in $h$, we first notice that $G(0, \kappa_0) = 0$. 
Writing $f(x) = \re^{-(\kappa_0 + \alpha)x}$ and applying Lemma \ref{l:Taylor-expansion}--\ref{i:taylor-1}, we subsequently compute 
\begin{align*}
G(h, \kappa_0 + \alpha) = {} & G(h, \kappa_0 + \alpha) - G(0, \kappa_0),\\[1ex]
= {}& -c \alpha - \kappa_0^2 + \sum_{k=1}^{+\infty} a_k \frac{f(kh) - 2 f(0) + f(-kh)}{k^2 h^2},\\
= {}& (2\kappa_0 - c)\alpha + \cO(h + \alpha^2).
\end{align*}
Since $\partial_\kappa G(0, \kappa_0) = 2\kappa_0 - c < 0$, we may apply the implicit function theorem to construct the claimed function $h\mapsto \kappa(h)$. Incidentally, it expands as $\kappa(h) = \kappa_0 + h^2 \frac{{\kappa_0}^4}{12(c-2\kappa_0)}(\sum_{k\geq 1} a_k k^2) + \cO(h^4)$.
\end{proof}

\begin{lemma}
\label{l:residual-bounds}
For any $\delta>0$, there exists $h_0>0$ such that for all $h\in (0, h_0)$,
\begin{equation*}
\norm{\cR}_{L^2\cap L^\infty} \leq \delta.
\end{equation*}
\end{lemma}
\begin{proof}
Our starting goal is to show that $\cR$ is exponentially localized. Restricting to $x\geq 1$, we see as stated previously that $\omega(x)^{-1} T_{2,g}(0, \phi^\infty(x)) = \cO(e^{(\theta - 2\kappa(h))x}) = \cO(e^{-\delta x})$. Continuing with $x\geq 1$, we turn to the infinite range linear term, and rely on Lemma \ref{l:definition-kappa} to express it as 
\begin{equation*}
\omega^{-1} \cA_h(0) \phi^\infty = \omega^{-1} \cA_h(0) (w_0 + \mathbf{1}_- (1 + \re^{-\kappa \cdot})).
\end{equation*}
Since $w_0$ is constructed from the solution to an ODE that converges to a hyperbolic equilibrium, it has exponential behavior at infinity, in the sense that
\begin{equation*}
\sup_{\absolute{s}\leq kh} \Absolute{\frac{w_0^{(3)}(x + s)}{e^{-(\kappa_0 + \delta) x}}} \leq C e^{(\kappa_0 + \delta) kh}.
\end{equation*}
We can obviously obtain a similar bound for $(1 + e^{-\kappa \cdot}) \mathbf{1}_- $, which allows us to apply Lemma \ref{l:Taylor-expansion}--\ref{i:taylor-2} and conclude that for all $x\geq 1$ we have
\begin{equation*}
\absolute{\omega^{-1} \cA_h(0) \phi^\infty}(x) \leq C e^{\theta x} e^{-(\kappa_0 + \delta) x} \leq C e^{-\frac{\delta}{2}x}.
\end{equation*}

We use a similar approach when $x\leq -1$, and remark that $\cR$ then simplifies to 
\begin{equation*}
\cR = \Delta_{a,h} (\mathbf{1}_- + w_0 + \mathbf{1}_+ e^{-\kappa(h)\cdot}) + c w_0' + T_{1,g}(1,w_0).
\end{equation*}
The infinite range term is controlled using Lemma \ref{l:Taylor-expansion}--\ref{i:taylor-2}, while the other two terms have the same exponential decay as $w_0$.

In particular, we conclude that $\absolute{\cR} \leq \re^{-\eta \absolute{x}}$ for some $\eta>0$, ensuring that $\cR$ belongs to $L^2\cap L^\infty$. To prove its smallness, we remark that $\cR$ converges to $0$ when $h\to 0$. Indeed, the smoothness of $\phi^\infty$ allows us to take the limit $h\to 0$ in the expression of $\cR$, leading to 
\begin{equation*}
\lim_{h\to 0} \cR = \omega^{-1} \bigg( \cA(0)\phi_0 + T_{2,g}(0,\phi_0)\bigg) = 0.
\end{equation*}
As a consequence, we can restrict to a sufficiently small value of $h_0$ to obtain the claim.
\end{proof}

\subsection{Finite differences}
In the following, we will use more finite difference operators. We denote the upwind and downwind first order differences as
\begin{equation*}
\partial_{kh}^+ v \deq \frac{v(\cdot + kh) - v}{kh},
\hspace{4em}
\partial_{kh}^- v \deq \frac{v - v(\cdot - kh)}{kh},
\end{equation*}
and their centered counterpart as
\begin{equation*}
\partial_{kh}^0 v \deq \frac{v(\cdot + kh) - v(\cdot - kh)}{2kh}.
\end{equation*}
We combine them to define the infinite range transport term
\begin{equation*}
\partial_{a,h}^0 v \deq \sum_{k=1}^{+\infty} a_k \partial_{kh}^0 v.
\end{equation*}
Let us further define the mean operator $M_{kh}^0 v := \frac{v(\cdot + kh) + v(\cdot - kh)}{2}$ together with its infinite range counterpart
\begin{equation*}
M_{a,h}^0 v \deq \sum_{k=1}^{+\infty} a_k \frac{v(\cdot + kh) + v(\cdot - kh)}{2}.
\end{equation*}
We start by summarizing several long but useful elementary computations.
\begin{lemma}
\label{l:finite-difference-computations}
Let $f\in \buc{3}(\RR, \RR)$ and $v\in L^2$. Then the following statements hold for all $0 < h \le 1$.
\begin{enumerate}
\item \label{i:finite-diff-1} (unbalanced mean) 
\begin{align*}
(\partial_{kh}^+ f) v(\cdot + kh) + (\partial_{kh}^- f) v(\cdot - kh) = {}& 2 f' M_{kh}^0 v\\
&{} + \bigg(T_{2,f}(\cdot, kh) v(\cdot + kh) - T_{2,f}(\cdot, -kh)v(\cdot - kh)\bigg)
\end{align*}
and the final term satisfies the remainder bound
\begin{equation}
\label{e:bound-remainder-2}
\frac{1}{kh} \Norm{T_{2,f}(\cdot, kh) v(\cdot + kh) - T_{2,f}(\cdot, -kh) v(\cdot - kh)}_{L^2} \leq kh \norm{f}_{W^{2,\infty}} \norm{v}_{L^2};
\end{equation}
\item \label{i:finite-diff-2} (unbalanced difference) 
\begin{align*}
\frac{1}{kh}\bigg((\partial_{kh}^+f) & v(\cdot + kh) - (\partial_{kh}^-f) v(\cdot - kh)\bigg)\\
= {}& 2 f' \partial_{kh}^0v  + f'' M_{kh}^0 v + \frac{1}{k^2 h^2}\bigg(v(\cdot + kh) T_{3,f}(\cdot, kh) + v(\cdot - kh) T_{3,f}(\cdot, -kh)\bigg)
\end{align*}
and the last term satisfies the remainder bound
\begin{equation}
\label{e:bound-remainder-1}
\frac{1}{k^2 h^2}\Norm{v(\cdot + kh) T_{3,f}(\cdot, kh) + v(\cdot - kh) T_{3,f}(\cdot, -kh)}_{L^2} \leq kh \norm{f}_{W^{3,\infty}} \norm{v}_{L^2};
\end{equation}
\item \label{i:finite-diff-3} (integration by parts)
\begin{align*}
\scalp{f\partial_{kh}^0 v, v} = {}& -\frac{1}{2}\scalp{f' M_{kh}^0 v, v} - \frac{1}{4 kh}\Scalp{T_{2,f}(\cdot, kh) v(\cdot + kh) - T_{2,f}(\cdot, -kh) v(\cdot - kh),v};
\end{align*}
\item \label{i:finite-diff-4}(upper semicontinuity)
\begin{equation*}
\norm{\partial_{kh}^0 v}_{L^2} \leq \norm{v'}_{L^2}.
\end{equation*}
\end{enumerate}
\end{lemma}
\begin{proof}
To obtain \ref{i:finite-diff-1}--\ref{i:finite-diff-2}, we simply Taylor expand $f$ up to order two and three respectively. The bounds are direct using the integral formulation of the Taylor remainder. For example, the identity
\begin{equation*}
T_{3, f}(x, kh) = \int_x^{x+kh} f^{'''}(y) \frac{(x + kh - y)^2}{2} \dd y
\end{equation*}
leads to the bound
\begin{align*}
(kh)^{-2} \norm{v(\cdot + kh)T_{3,f}(\cdot, kh)}_{L^2} \leq {}& (kh)^{-2} \cdot \frac{1}{6} (kh)^3 \norm{f}_{W^{3,\infty}} \cdot \norm{v}_{L^2} \\
\leq {}& \frac{1}{2} kh \norm{f}_{W^{3,\infty}} \norm{v}_{L^2}.
\end{align*}

Turning to \ref{i:finite-diff-3}, we recall the summation by parts
identity $\scalp{\partial_{kh}^\pm u, v} = - \scalp{u, \partial_{kh}^\mp v}$
together with the Leibniz rule $\partial_{kh}^\pm(uv) = (\partial_{kh}^\pm u) v + u(\cdot \pm kh) (\partial_{kh}^\pm v)$. We first compute 
\begin{align}
\nonumber \scalp{f\partial_{kh}^0 v, v} = {} & - \scalp{v, \partial_{kh}^0 (f v)},\\ 
\label{e:IBP-first-order} = {} & - \scalp{v, f\partial_{kh}^0 v} - \frac{1}{2}\scalp{v, (\partial_{kh}^+ f) v(\cdot + kh) + (\partial_{kh}^- f) v(\cdot - kh)}
\end{align}
and recognize the expression from \ref{i:finite-diff-1}. Rearranging, we obtain the claimed equality.

Finally, we set out to prove \ref{i:finite-diff-4}. We use Plancherel's theorem and the inequality $\absolute{\sin(x)}\leq \absolute{x}$ to obtain
\begin{equation*}
\norm{\partial_{kh}^0 v}_{L^2} = \frac{1}{kh}\norm{\xi \mapsto \ri \sin(kh\xi) \hat{v}(\xi)}_{L^2} \leq \norm{\xi \mapsto \ri \xi \hat{v}(\xi)}_{L^2} = \norm{v'}_{L^2},
\end{equation*}
which is the claimed bound.
\end{proof}

\section{Linear theory}
\label{s:linear-theory}

Our goal in this section is to prove invertibility of the weighted linearization, stated as
\begin{proposition}
\label{p:Lh-is-invertible}
There exists $h_0>0$ such that for all $h\in (0, h_0)$ and all $f\in L^2(\RR)$, the resolvent problem 
\begin{equation*}
\cL_h v = f
\end{equation*}
admits a unique solution $v\in H^1(\RR)$. Furthermore, there exists a positive constant $C$ such that for all $h$ and $f$,
\begin{equation*}
\norm{\cL_h^{-1} f}_{H^1} \leq C \norm{f}_{L^2}.
\end{equation*}
\end{proposition}

To prove this statement, we start by studying the space-continuous setting, and then transfer properties from the linear operator to $\cL_h$. The transfer procedure will take most of our attention. It will require us to first simplify the expression of $\cL_h$, and then to obtain resolvent bounds that are uniform in $h$.

\subsection{The continuous case}
Our starting point is the study of the corresponding second order differential operator
\begin{equation*}
\cL v \deq \omega^{-1} \cA(\phi^\infty) (\omega v) \deq \omega^{-1} \bigg(\partial_{xx} + c\partial_x + g'(\phi^\infty)\bigg) (\omega v),
\end{equation*}
together with its adjoint $\cL^\ad$.
\begin{lemma}
\label{l:L-is-invertible}
For all $f\in L^2(\RR)$, there exists a unique solution $v\in H^2(\RR)$ to 
\begin{equation*}
\cL v = f.
\end{equation*}
Furthermore, there exists a positive constant $C$ such that for all $f$
\begin{equation*}
\norm{\cL^{-1}f}_{H^2} \leq C \norm{f}_{L^2}.
\end{equation*}
The same holds for the adjoint: there exists $C>0$ such that for every $f\in L^2(\RR)$, the adjoint eigenproblem has a unique solution in $H^2(\RR)$, which satisfies
\begin{equation*}
\norm{(\cL^\ad)^{-1}f}_{H^2} \leq C \norm{f}_{L^2}.
\end{equation*}
\end{lemma}
\begin{proof}
Let us start with a detailed proof of the statements for $\cL$. We actually show that there exists $\delta>0$ such that the spectrum of $\cL$ is a subset of the stable sector
\begin{equation*}
S(\delta) \deq \Set{\lambda \in \CC : \Real{\lambda} < -\delta - \delta \absolute{\Imag{\lambda}}}.
\end{equation*}
This implies that $0$ does not belong to the spectrum of $\cL$, which is precisely the statement. To do so, we even show that $S(\delta)$ contains the numerical range of $\cL$, which is the subset of $\CC$ defined as 
\begin{equation*}
\Set{\Scalp{\cL u, u} : \norm{u}_{L^2} = 1},
\end{equation*}
which contains the spectrum of $\cL$, see \emph{e.g.} \cite[Lemma 4.1.9]{Kapitula-Promislow-13}.
We first compute that 
\begin{equation}
\label{e:expression-L}
\cL u = \partial_{xx} u + \left(c + 2\frac{\omega'}{\omega}\right) \partial_x u + \left(\frac{\omega''}{\omega} + c\frac{\omega'}{\omega} + g'(\phi^\infty)\right) u.
\end{equation}
Testing this equality against $u$, integrating by part the second order derivative, and recalling for any differentiable $h:\RR\to \RR$ the identity $\Real{\scalp{h\partial_x u, u}} = -\frac{1}{2}\scalp{h' u, u}$, we compute
\begin{equation*}
\Real{\Scalp{\cL u, u}} = -\norm{\partial_x u}_{L^2}^2 + \Scalp{\left(\left(\frac{\omega'}{\omega}\right)^2 + c\frac{\omega'}{\omega} + g'(\phi^\infty)\right)u, u}.
\end{equation*}

We now proceed by replacing $\omega$ with $\tilde{\omega}$; see \eqref{e:weight} and the discussion that follows. For any $x\in (-1,1)$, we can use the fact that $g$ is concave to estimate
\begin{align*}
\left(\frac{\tilde{\omega}'}{\tilde{\omega}}\right)^2 + {} & c\frac{\tilde{\omega}'}{\tilde{\omega}} + g'(\phi^\infty) \\
= {} & \theta^2 - c\theta + g'(0) + \frac{\theta^2}{4} \left((x+1)^2 - 4\right) + c\frac{\theta}{2} \left(x+1 - 2\right) + \left(g'(\phi^\infty) - g'(0)\right), \\
\leq {} & \theta^2 - c\theta + g'(0).
\end{align*}
Thus, for all $x \in \RR$ we obtain 
\begin{equation*}
\left(\frac{\tilde{\omega}'}{\tilde{\omega}}\right)^2 + c\frac{\tilde{\omega}'}{\tilde{\omega}} + g'(\phi^\infty) \leq \max\left(g'(1), \theta^2 - c\theta + g'(0)\right).
\end{equation*}
We remark that the right hand side is negative due to assumption \ref{a:g} and the choice of $\theta$. As a consequence, we can ensure by taking $\varepsilon$ sufficiently small in the definition of $\omega$ that there exists $\delta>0$ such that
\begin{equation*}
\left(\frac{\omega'}{\omega}\right)^2 + c\frac{\omega'}{\omega} + g'(\phi^\infty) \leq - \delta.
\end{equation*}
In particular, we see that $\Real{\scalp{\cL u, u}} \leq - \norm{\partial_x u}_{L^2}^2 - \delta \norm{u}_{L^2}^2$. Turning to the imaginary part, we simply estimate
\begin{equation*}
\absolute{\Imag{\scalp{\cL u, u}}} = \absolute{\Imag{\scalp{(c + 2\omega^{-1} \omega') \partial_x u, u}}} \leq C \norm{\partial_x u}_{L^2} \norm{u}_{L^2}. 
\end{equation*}
As a result, we obtain the claimed bound from
\begin{equation*}
\Real{\scalp{\cL u, u}} \leq -\norm{\partial_x u}_{L^2} - \delta \leq - \frac{1}{C} \Absolute{\Imag{\scalp{\cL u, u}}} - \delta.
\end{equation*}
We conclude the proof by considering $\cL^\ad$. By standard theory, its spectrum is the conjugate of the spectrum of $\cL$. In particular, $0$ does not belong to the spectrum of $\cL^\ad$, and the proof is complete.
\end{proof}

\subsection{Reformulation of the discrete weighted operator}
Before being able to transfer invertibility, we first express $\cL_h$ in a fashion that is comparable to \eqref{e:expression-L}. For this, we need to extract the weight out of the finite differences and to isolate the relevant terms.
\begin{proposition}
\label{p:Delta-h-decomposition}
For all $v\in L^2(\RR)$ we have the expansion
\begin{equation*}
\omega^{-1}\Delta_{a,h}(\omega v) = \Delta_{a,h} v + 2\frac{\omega'}{\omega} \partial_{a,h}^0 v + \frac{\omega''}{\omega} M_{a,h}^0 v + Rv,
\end{equation*}
in which the residual satisfies the bound
\begin{equation}
\label{e:bound-residual}
\norm{Rv}_{L^2} \leq C h \norm{v}_{L^2}.
\end{equation}
\end{proposition}
\begin{proof}
We start with the decomposition
\begin{equation}
\label{e:decomposition-Delta-1}
\Delta_h(\omega v) = \omega \Delta_{a,h} v + \sum_{k=1}^{+\infty}a_k \frac{1}{kh}\bigg(v(\cdot + kh) \partial_{kh}^+\omega - v(\cdot - kh) \partial_{kh}^- \omega \bigg),
\end{equation}
and apply Lemma \ref{l:finite-difference-computations}--\ref{i:finite-diff-2} to obtain the stated expansion, with the associated remainder
\begin{align*}
R = \omega^{-1}\sum_{k=1}^{+\infty} a_k \frac{1}{k^2h^2}\bigg(v(\cdot + kh) T_{3,\omega}(\cdot, kh) + v(\cdot - kh) T_{3,\omega}(\cdot, -kh)\bigg).
\end{align*}
The bound \eqref{e:bound-remainder-1} is almost sufficient 
for our purposes here. We refine it by incorporating the $\omega^{-1}$ into it as was done in Lemma \ref{l:Taylor-expansion}. The exponential behavior of $\omega$ at infinity ensures that
\begin{equation*}
\sup_{\absolute{y-x}\leq kh} \frac{\omega'''(y)}{\omega(x)} \leq \Norm{\frac{\omega'''}{\omega}}_{L^\infty} \ \sup_{\absolute{y-x}\leq kh} \frac{\omega(y)}{\omega(x)} \leq C e^{(\kappa_0 + \delta) kh},
\end{equation*}
which in turn implies
\begin{align*}
\norm{\omega^{-1} T_{3,\omega}(\cdot, kh)}_{L^\infty} \leq C k^3h^3 e^{(\kappa_0 + \delta) kh}. 
\end{align*}
Such a sequence grows slowly enough when $h$ is taken to be sufficiently small that the sum against $(a_k)_k$ converges. This concludes the proof.
\end{proof}

As a consequence of the result above, we are able to express $\cL_h$ in the form
\begin{equation}
\label{e:cL-h}
\cL_h v = \Delta_{a,h} v + c v' + 2\frac{\omega'}{\omega} \partial_{a,h}^0 v + \left(g'(\phi^\infty) + c\frac{\omega'}{\omega} + \frac{\omega''}{\omega} M_{a,h}^0 \right)v+ R v,
\end{equation}
where the remainder operator $R$ satisfies \eqref{e:bound-residual}. It is now convenient to compute the adjoint operator $\cL_h^\ad$ using this formulation. Computations similar to \eqref{e:IBP-first-order} and those of Proposition \ref{p:Delta-h-decomposition} lead to the expression
\begin{equation}
\label{e:cL-h-ad}
\cL_h^\ad v = \Delta_{a,h} v - c v' - 2\frac{\omega'}{\omega} \partial_{a,h}^0 v + \left(g'(\phi^\infty) + c\frac{\omega'}{\omega} + \frac{\omega''}{\omega} M_{a,h}^0 - 2\left(\frac{\omega'}{\omega}\right)'\right) v + R^\ad v,
\end{equation}
with a remainder operator $R^\ad$ that satisfies the same bound \eqref{e:bound-residual}. Let us emphasize that $R^\ad$ does not represent the adjoint of $R$.

\subsection{Resolvent bounds}

This part mostly follow \cite{Hupkes-Jukic-24}. Let us introduce
\begin{equation*}
\Lambda(h) \deq \inf \Set{\norm{\cL_h v}_{L^2} : v \in H^1(\RR), \norm{v}_{H^1} = 1},
\end{equation*}
together with 
\begin{equation*}
\Lambda \deq \liminf_{h\downarrow 0} \Lambda(h).
\end{equation*}
In the next subsection, we will actually need similar definitions involving the adjoint operator of $\cL_h$. To simplify the discussion, we refrain to provide these here, and will explain at a later stage how to adapt the content that follows.

Our goal is to show that $\Lambda>0$. We establish a bound from below by constructing a weakly-convergent sequence.

\begin{lemma}
\label{l:construct-sequence}
There exists a sequence $(h_n, v_n, f_n)_{n\geq 0}$ of $(0, \delta_0) \times H^1 \times L^2$, together with $v_* \in H^1$ and $f_* \in L^2$, such that the following holds.
\begin{enumerate}
\item \label{i:construct-seq-1} For all $n\geq 0$, $\norm{v_n}_{H^1} = 1$ and 
\begin{equation}
\label{e:eigenproblem-seqence}
\cL_{h_n} v_n = f_n.
\end{equation}
When $n\to +\infty$, we have 
\begin{align*}
h_n \longrightarrow 0, \\
\norm{f_n}_{L^2} \longrightarrow \Lambda.
\end{align*}
\item \label{i:construct-seq-2} When $n\to +\infty$, the following weak convergences hold:
\begin{align*}
v_n \rightharpoonup {}& v_* & \text{in } H^1, \\
f_n \rightharpoonup {}& f_* & \text{in } L^2, \\
\cL_{h_n} v_n \rightharpoonup {}& \cL v_* & \text{in } L^2.
\end{align*}
For any $m>0$, when $n\to +\infty$, the following strong convergence holds: 
\begin{equation*}
\norm{v_n}_{L^2(-m,m)} \to \norm{v_*}_{L^2(-m,m)}.
\end{equation*}
\item \label{i:construct-seq-3} There exists a positive constant $C_1$ such that 
\begin{equation*}
\Lambda \geq C_1 \norm{v_*}_{H^2}.
\end{equation*}
\end{enumerate}
\end{lemma}
\begin{proof}
The existence of a sequence $(h_n, v_n, f_n)$ satisfying \ref{i:construct-seq-1} is a direct consequence of the definition of $\Lambda$. Since $\cL_{h_n}$ has real coefficients, we can ensure that $v_n$ and $f_n$ have real values.

We now prove \ref{i:construct-seq-2}. Using the first point, together with equation \eqref{e:eigenproblem-seqence}, we successively get that the sequences $(v_n, f_n)_{n\geq 0}$ and $(v_n, f_n, \cL_{h_n} v_n)_{n\geq 0}$ are bounded. Extracting a sub-sequence, we obtain the weak convergence of the latter sequence in $H^1 \times L^2 \times L^2$ to some element $(v_*, f_*, y_*)$. 
We show that $y_* = \cL v_*$ by testing against a smooth and compactly supported function $\xi$:
\begin{align*}
\lim_{n\to +\infty} \scalp{\cL_{h_n} v_n, \xi} = {} & \lim_{n\to +\infty} \scalp{v_n, \cL_{h_n}^\ad \xi}, \\ 
= {} & \scalp{v_*, \cL^\ad\xi},\\
= {} & \scalp{\cL v_*, \xi},
\end{align*}
the second equality coming from the fact that $\cL_{h_n}^\ad \xi$ strongly converges to $\cL^\ad \xi$. Uniqueness of the weak limit concludes this point. It remains to show the strong $L^2$-convergence on compact sets. This follows from the weak $H^1$-convergence and the fact that the injection $H^1(-m, m) \hookrightarrow L^2(-m,m)$ is compact.

Turning to \ref{i:construct-seq-3}, we pass to the weak limit in \eqref{e:eigenproblem-seqence} to obtain $\cL v_* = f_*$. Using the invertibility of $\cL$, the weak lower semicontinuity of the $L^2$-norm and the convergence of $\norm{f_n}_{L^2}$, we successively obtain
\begin{equation*}
C \norm{v_*}_{H^2} \leq \norm{f_*}_{L^2} \leq \liminf_{n\to +\infty} \norm{f_n}_{L^2} = \Lambda
\end{equation*}
for a positive constant $C$, thus concluding the proof.
\end{proof}

To conclude that $\Lambda \neq 0$, we now rely on the specific structure of the problem to bound $\norm{v_*}_{L^2}$ from below. We test the eigenvalue equation against $v'$ and $v$ respectively, and handle the different terms that appear in these expressions. 
\begin{lemma}
\label{l:no-oscillations}
There exists positive constants $C$ and $h_0$ such that for all $h\in (0, h_0)$ and for all $(v,f)\in H^1\times L^2$ that satisfy $\cL_h v = f$, 
\begin{equation*}
\norm{v'}_{L^2}^2 \leq C \norm{f}_{L^2}^2 + C \norm{v}_{L^2}^2.
\end{equation*}
\end{lemma}
\begin{proof}
We test the eigenproblem $\cL_h v = f$ against $v'$, and use \eqref{e:cL-h} to obtain
\begin{equation*}
c\norm{v'}_{L^2}^2 \leq \norm{f}_{L^2} \norm{v'}_{L^2} + \Absolute{\Scalp{\Delta_{a,h}v, v'}} + 2 \Norm{\frac{\omega'}{\omega}}_{L^\infty} \norm{\partial_{a,h}^0 v}_{L^2} \norm{v'}_{L^2} + C \norm{v}_{L^2} \norm{v'}_{L^2}. 
\end{equation*}
We now successively inspect the different terms, starting with the second order difference. We claim that $\scalp{\Delta_{a,h}v, v'} = 0$: indeed for every $k\geq 1$,  integration and summation by parts gives
\begin{equation*}
\scalp{\Delta_{kh} v, v'} = - \scalp{v', \Delta_{kh} v} = 0. 
\end{equation*}
Proceeding to the next term, we remark that the definition of $\omega$ ensures that $\norm{\omega^{-1}\omega'}_{L^\infty} \leq \theta + \varepsilon$; see Lemma \ref{l:construct-omega}. In combination with Lemma \ref{l:finite-difference-computations}--\ref{i:finite-diff-4}, this leads to 
\begin{equation*}
(c - 2\theta - 2\varepsilon) \norm{v'}_{L^2}^2 \leq \norm{f}_{L^2} \norm{v'}_{L^2} + C \norm{v}_{L^2} \norm{v'}_{L^2}.
\end{equation*}
At this stage, it is useful to take $\varepsilon$ so small that $c-2\theta-2\varepsilon > 0$, which is possible due to Lemma \ref{l:construct-theta}. We now apply Young's inequality $xy \leq \frac{\alpha}{3}x^2 + \frac{3}{4\alpha}y^2$ twice with $\alpha = c-2\theta-2\varepsilon$ and $x = \norm{v'}_{L^2}$ to arrive at the claimed estimate.
\end{proof}

\begin{proposition}
\label{p:reverse-inequality}
There exist positive constants $m$, $C_2$ and $C_3$ such that 
\begin{equation*}
\norm{v_*}_{L^2(-m,m)}^2 \geq C_2 - C_3 \Lambda^2.
\end{equation*}
\end{proposition}
\begin{proof}
We test the eigenproblem against $v$, once again using \eqref{e:cL-h} to obtain
\begin{equation*}
\scalp{f, v} = \scalp{\Delta_{a,h} v, v} + c\scalp{v', v} + 2\Scalp{\frac{\omega'}{\omega}\partial_{a,h}^0 v, v} + \Scalp{\left(g'(\phi^\infty) + c\frac{\omega'}{\omega} + \frac{\omega''}{\omega}M_{a,h}^0 + R\right)v, v}.
\end{equation*}
We see that for all sufficiently small $h$: 
\begin{equation*}
\scalp{\Delta_{a,h} v, v} = 2 \sum_{k=1}^{+\infty} a_k \frac{\cos(kh) - 1}{k^2 h^2}\scalp{\hat{v}, \hat{v}} \leq 0.
\end{equation*}
Indeed, dominated convergence and assumption \ref{a:a} ensures that $2\sum_{k=1}^{+\infty} a_k \frac{\cos(kh) - 1}{k^2 h^2} \to -1$ when $h\to 0$.
Further using the identity $\scalp{v, v'} = 0$ (which is easily obtained by integrating by parts once), and applying Lemma \ref{l:finite-difference-computations}--\ref{i:finite-diff-3} with $f=2\frac{\omega'}{\omega}$, we obtain
\begin{equation}
\label{e:reverse-inequality-3}
\scalp{f, v} \leq \Scalp{\left(g'(\phi^\infty) + c\frac{\omega'}{\omega} + \left(\frac{\omega'}{\omega}\right)^2 M_{a,h}^0 + \tilde{R}\right)v, v},
\end{equation}
where the remainder satisfies $\norm{\tilde{R} v}_{L^2} \leq C h \norm{v}_{L^2}$ thanks to \eqref{e:decomposition-Delta-1} and \eqref{e:bound-remainder-2}.
At this stage, it is worth mentioning that while $\Scalp{\frac{\omega'}{\omega} v, v}$ has a sign, this is not the case for $\Scalp{\left(\frac{\omega'}{\omega}\right)^2 M^0_{a,h}v, v}$. Thus, extra care should be taken with the latter.

Using Assumption \ref{a:g}, the definition of $\omega$ and \eqref{e:stable-spectrum}, we see that both
\begin{equation*}
\lim_{-\infty} g'(\phi^\infty) + c\frac{\omega'}{\omega} + \left(\frac{\omega'}{\omega}\right)^2 + Ch = g'(1) + Ch
\end{equation*}
and 
\begin{equation*}
\lim_{+\infty} g'(\phi^\infty) + c\frac{\omega'}{\omega} + \left(\frac{\omega'}{\omega}\right)^2 + Ch = g'(0) - c\theta + \theta^2 + Ch
\end{equation*}
are negative, possibly after restricting to smaller values of $h$. 
For any $m \in \RR$, we use the second limit to compute 
\begin{align*}
\scalp{f,v}_{L^2(m,+\infty)} \leq {}& \bigg( g'(0) - c\theta + \theta^2 + Ch + \hspace{-1ex} \smallo_{\makebox[0pt]{\phantom{I}} m \to \mathrlap{+\infty}}(1) \bigg) \norm{v}^2_{L^2(m,+\infty)} \\
&{} + \Scalp{\left(\frac{\omega'}{\omega}\right)^2 M_{a,h}^0 v - \theta^2 v, v}_{L^2(m,+\infty)}.
\end{align*}
Since $\frac{\omega'}{\omega} = -\theta$ for large enough $m$, a simple Cauchy-Schwartz inequality allows to bound the second term as
\begin{equation*}
\theta^2 \left(\Scalp{M_{a,h}^0 v, v}_{L^2(m,+\infty)} - \norm{v}^2_{L^2(m,+\infty)} \right)\leq 0.
\end{equation*}
Setting $m$ large enough in the above discussion, and using very similar computations on the set $(-\infty, -m)$, we conclude from \eqref{e:reverse-inequality-3} that there exist a positive constant $\eta$ such that
\begin{equation}
\label{e:reverse-inequality-1}
-\norm{f}_{L^2} \norm{v}_{L^2} \leq \scalp{f, v} \leq C \norm{v}_{L^2(-m,m)}^2 - \eta \norm{v}_{L^2}^2.
\end{equation}

We now apply \ref{e:reverse-inequality-1} to the sequence constructed in Lemma \ref{l:construct-sequence}. Invoking Young's inequality, we obtain
\begin{align}
\nonumber
C \norm{v_n}_{L^2(-m,m)}^2 \geq {}& \eta \norm{v_n}_{L^2}^2 - \norm{f_n}_{L^2} \norm{v_n}_{L^2},\\
\label{e:reverse-inequality-2}
\geq {}& \frac{\eta}{2} \norm{v_n}_{L^2}^2 - \frac{1}{2\eta} \norm{f_n}_{L^2}^2.
\end{align}
We then add to \eqref{e:reverse-inequality-2} a $-\beta$ multiple of the bound from Lemma \ref{l:no-oscillations}, choosing $\beta = \frac{\eta}{2(C+1)}$ as the solution of $\frac{\eta}{2} - \beta C = \beta$ to obtain
\begin{align*}
C \norm{v_n}_{L^2(-m,m)}^2 \geq {}& \left(\frac{\eta}{2} - \beta C\right) \norm{v_n}_{L^2}^2 + \beta\norm{v_n'}_{L^2} - \left(\frac{1}{2\eta} + C\beta\right) \norm{f_n}_{L^2}^2, \\
\geq {}& \beta - C \norm{f_n}_{L^2}^2.
\end{align*}
To conclude the proof, we recall that $\norm{f_n} \leq \Lambda$ was obtained during the proof of Lemma \ref{l:construct-sequence}--\ref{i:construct-seq-3}, and use the strong convergence of $v_n$ towards $v_*$ in $L^2(-m,m)$ when $n\to +\infty$.
\end{proof}

We can finally combine Lemma \ref{l:construct-sequence}--\ref{i:construct-seq-3} and Proposition \ref{p:reverse-inequality}, to obtain
\begin{equation*}
\Lambda^2 \geq C_1 \norm{v_*}_{L^2}^2 \geq C_1 \norm{v_*}_{L^2(-m,m)}^2 \geq C_1 C_2 - C_1 C_3 \Lambda^2.
\end{equation*}
Rearranging, we find
\begin{equation*}
\Lambda^2 \geq \frac{C_1 C_2}{1 + C_1 C_3} > 0,
\end{equation*}
establishing that indeed $\Lambda>0$. 

To conclude this section, we explain how the material above can be adapted to the adjoint problem. Proceeding from the representation \eqref{e:cL-h-ad}, we introduce the quantities
\begin{equation*}
\Lambda^\ad(h) \deq \inf \Set{\norm{\cL_h^\ad v}_{L^2} : v \in H^1(\RR), \norm{v}_{H^1} = 1},
\end{equation*}
and
\begin{equation*}
\Lambda^\ad \deq \liminf_{h\downarrow 0} \Lambda^\ad(h).
\end{equation*}
which are analogous to those analyzed above for $\mathcal{L}_h$.
It is then straightforward to adapt Lemmas \ref{l:construct-sequence}, Lemma \ref{l:no-oscillations} and Proposition \ref{p:reverse-inequality} to the adjoint case. From there, $\Lambda^\ad>0$ follows as above.
\subsection{The discrete case: proof of Proposition \ref{p:Lh-is-invertible}}

We are finally ready to prove the main statement of section \ref{s:linear-theory}. We follow \cite{Bates-Chen-Chmaj-03}.

\begin{proof}[Proof of Proposition \ref{p:Lh-is-invertible}]
Since $\Lambda > 0$, there exists $h_0>0$ such that $\Lambda(h)>0$ for all $h\in (0, h_0)$. For such $h$, it follows that $\cL_h$ has a bounded inverse, when seen as linear operator from $H^1$ to its range $\cR \subset L^2$. In particular, $\cR$ is closed in $L^2$, and we assume by contradiction that $\cR \neq L^2$. 

Then there exists a nonzero $w\in L^2$ such that for all $v\in H^1$, 
\begin{equation*}
0 = \scalp{\cL_h v, w} = \scalp{v, \cL_h^\ad w},
\end{equation*}
thus ensuring $\cL_h^\ad w = 0$. This is a contradiction with the fact that $\Lambda^\ad(h) > 0$. The bound $\norm{\cL_h^{-1} f}_{H^1} \leq \frac{1}{\Lambda(h)} \norm{f}_{L^2}$ comes from the definition of $\Lambda(h)$. 
\end{proof}

\begin{remark}
\label{r:h-0-value}
In the previous proof, control of $h_0$ is nearly intractable. It would require estimates involving the behaviour of $\Lambda(h)$ as $h$ approaches $0$. This is out of reach for now.
\end{remark}

\section{Fix point argument}
We are now ready to prove our main result concerning the existence of travelling fronts. As a final preparation, we obtain an estimate for the quadratic terms $\cQ$. 

\begin{lemma}
\label{l:nonlinear-bounds}
Let $p\in\Set{2,+\infty}$. There exists $C>0$ and $h_0>0$ such that for all $h\in (0, h_0)$, and all $v$ such that $\norm{v}_{L^\infty} \leq 1$,
\begin{equation*}
\norm{\cQ(v_1) - \cQ(v_2)}_{L^p} \leq \norm{v_1 - v_2}_{L^p} \left(\norm{v_1}_{L^\infty} + \norm{v_2}_{L^\infty}\right).
\end{equation*}
\end{lemma}
\begin{proof}
We use the Taylor expansion structure of $\cQ$. For every set of reals $a$, $b_1$ and $b_2$, we may use the integral formulation to compute
\begin{align*}
T_{2,g}(a,b_2) - T_{2,g}(a,b_1) = {}& \int_a^{a+b_2} \frac{g''(s)}{2!} (a + b_2 - s) \dd s - \int_a^{a+b_1} \frac{g''(s)}{2!} (a + b_1 - s) \dd s, \\
= {}& \int_{a+b_1}^{a+b_2} \frac{g''(s)}{2!} (a + b_2 - s) \dd s - \int_a^{a+b_1} \frac{g''(s)}{2!} (b_2 - b_1) \dd s.
\end{align*}
Both of these integrals can be bounded by $C \absolute{b_2-b_1}\left(\absolute{b_1} + \absolute{b_2}\right)$, which directly leads to the claimed estimate.
\end{proof}

\begin{proof}[Proof of Theorem \ref{t:local}]
We reformulate Lemma \ref{l:v-formulation}, and see that it is enough to construct a solution $v\in L^\infty(\RR)$ to the fix point equation
\begin{equation*}
\label{e:fix-point}
v = \cL_h^{-1}\left(\cR + \cQ(v)\right).
\end{equation*}
Thanks to Lemma \ref{l:residual-bounds}, Proposition \ref{p:Lh-is-invertible} and Lemma \ref{l:nonlinear-bounds}, the map $v\mapsto \cL_h^{-1}\left(\cR + \cQ(v)\right)$ is a contraction on the Banach space 
\begin{equation*}
\Set{v\in L^2 \cap L^\infty : \norm{v}_{L^\infty} \leq \delta}, \quad  \norm{\cdot}_{L^2} + \norm{\cdot}_{L^\infty}
\end{equation*}
provided that $\delta$ and $h_0$ are chosen small enough. Applying the Banach fix point theorem concludes the proof.
\end{proof}

\bibliographystyle{alphaabbr}
\bibliography{discr-mono.bib}

\end{document}